\documentclass[12pt]{amsart}
\usepackage{amsmath,amssymb}

\newtheorem{thm}{Theorem}

\newtheorem{lemma}[thm]{Lemma}
\newtheorem{cor}[thm]{Corollary}

\newcommand{\E}{{\mathbb E}}

\newcommand{\normo}[1]{{\left\|#1\right\|}}

\newcommand{\snormo}[1]{{\mathopen\|#1\mathclose\|}}

\newcommand{\measure}{\text{measure}}

\begin{document}
\title[R.i. norms of symmetric sequence norms]
{Rearrangement Invariant Norms of Symmetric Sequence Norms
of Independent Sequences of Random Variables
}

\author{Stephen Montgomery-Smith}
\makeatletter
\address{Department of Mathematics\\
University of Missouri\\
Columbia, MO 65211}
\email{stephen@math.missouri.edu}
\urladdr{http://www.math.missouri.edu/\~{}stephen}
\thanks{The author was 
partially supported
by NSF grant DMS 9870026, and a grant from the Research Office of the
University of Missouri.}
\keywords{Independent random variables, symmetric sequence space,
rearrangement invariant space}
\subjclass{60G50, 46B45, 46E30}

\begin{abstract}
\noindent
Let $X_1$, $X_2,\dots,$\ $X_n$ be a sequence of independent random
variables, let $M$ be a rearrangement invariant space on the underlying
probability space, and let $N$ be a symmetric sequence space.  
This paper
gives an approximate formula for the quantity
$\snormo{\snormo{(X_i)}_N}_M$ whenever $L_q$ embeds into $M$ for some
$1 \le q < \infty$.
This extends work of Johnson and Schechtman who tackled the case when
$N = \ell_p$, and
recent work of Gordon, Litvak, Sch\"utt and Werner who obtained
similar results for Orlicz spaces.
\end{abstract}

\maketitle

\section{Introduction}

In a recent paper \cite{gordon et al}, Gordon, Litvak, Sch\"utt and Werner
considered the problem of computing $\E \normo{(a_i \xi_i)}_N$, where
$N$ is an Orlicz sequence space, $a_1$, $a_2,\dots,$\ $a_n$ are real numbers, 
and $\xi_1$, $\xi_2,\dots,$\ $\xi_n$ are identically distributed random 
variables.  
They were able to construct an Orlicz function $\Lambda$ such
that 
\[ 
   \E \normo{(a_i \xi_i)}_N \approx \snormo{(a_i)}_\Lambda ,
\]
where $A \approx B$ means that the ratio of $A/B$ is bounded below
and above by constants.  
Yehoram Gordon asked the question if a similar
formula could be found for general independent random variables --- 
not just scalar multiples of identically distributed random variables.
Mark Rudelson asked about the possibility of computing higher moments.

Seeking to answer these questions, the author found it to 
his psychological
advantage to consider the case of general rearrangement invariant spaces
rather than just Orlicz norms.
Let $M$ be a rearrangement invariant space on $[0,1]$, or equivalently
on any probability space, and let $N$ be a symmetric sequence space.
We will assume that all these vector spaces 
satisfy the triangle inequality.  
It seems quite likely that many
of the formulae will extend to at least some quasi-Banach situations,
but we do not explore this possibility here.  
We will normalize the
spaces so that $\snormo 1_M = \snormo{(1,0,\dots,0)}_N = 1$.
We seek to find an approximate formula for $\normo{\normo{(X_i)}_N}_M$.

The appropriate concept for describing our formula, the disjoint sum, 
or {\bf disjunctification}, has been
present in the literature for some time, for example, 
\cite{carothers-dilworth},
\cite{hitczenko-montgomery-smith},
\cite{johnson et al},
\cite{johnson-schechtman}. 
This is the function on $[0,n]$
that takes $t$ to $X^\#_{[t]+1}(t-[t])$, 
where $[t]$ denotes the integer part of $t$, and $X^\#_i$ 
denotes the non-increasing
rearrangment of $X_i$.  
We will write $Y$ for the
non-increasing rearrangement of this function, that is, 
$Y\colon[0,n] \to [0,\infty]$ is a non-increasing
function such that
\[ 
   \measure\{Y > t \} = \sum_{i=1}^n \Pr(X_i > t) . 
\]
The conjecture we will explore is
\begin{equation}
\label{e main}
   \snormo{\snormo{(X_i)}_N}_M \approx
   \snormo{Y|_{[0,1]}}_M + \snormo{(Y(i))}_N ,
\end{equation}
where the constants of approximation may be allowed to depend upon $M$.

Indeed for all the special cases hitherto considered in the literature,
this conjecture is true,
as long as $M$ is far away from $L_\infty$.
(It is clear that there must be some such restriction.  
For example,  
if $M = L_\infty$ then $\snormo{\snormo{(X_i)}_N}_M = 
\snormo{(\snormo{X_i}_\infty)}_N$, and so equation~(\ref{e main})
does not necessarily hold.)
First, Rosenthal's
inequality \cite{rosenthal} can be interpreted
as the truth of this statement in the cases that 
$N = \ell_1$ or $N = \ell_2$,
and $M = L_p$ for $1 \le p < \infty$.  
This was extended by
Carothers and Dilworth \cite{carothers-dilworth} to the case when
$M$ is a Lorentz space $L_{p,q}$, $1 \le p < \infty$, $1 \le q \le \infty$,
and then by Johnson and Schechtman \cite{johnson-schechtman} to the case
when $M$ is any rearrangement invariant space, even including the case
when $M$ is a quasi-Banach space, as long as $L_q$ embeds into $M$, and $M$ 
embeds into $L_r$, via the natural embeddings, for some $0<r\le q < \infty$.
It is not hard to extend this last result to also allow $N = \ell_p$ for any
$1 \le p < \infty$.

The main result of this paper is the following.

\begin{thm}
\label{t main}
Equation~(\ref{e main}) is true if there exists $1\le q < \infty$ such that 
$L_q$ embeds continuously via
the natural embedding into $M$.  
In that case, the constants of approximation
in equation~(\ref{e main})
depend only upon $q$ and the constant of embedding.
\end{thm}

In the last part of this paper,
we will describe how to use this to recover some of the results of 
\cite{gordon et al}.

The author would like to express his sincere appreciation to Mark Rudelson
for useful conversations, and 
for bringing this problem and reference \cite{gordon et al} to his attention,
and also to Joel Zinn for pointing out the reference \cite{marcus-zinn}.

\section{Proof of Main Theorem}

If $(x_i)$ is a sequence we will denote its non-increasing rearrangement
by $(x^*_i)$.  
If $f$ is a function or random variable, we will denote
its non-increasing rearrangement by $f^\#$.

\begin{lemma}
\label{l L_1 l_infty}
Equation~(\ref{e main}) is true if $M = L_1$ and $N = \ell_\infty$.
\end{lemma}

\begin{proof}
This follows because 
\begin{equation}
\label{e max in Pr}
    {\textstyle\frac12}\,\measure\{ Y|_{[0,1]} > t\}
    \le
    \Pr(\max_i X_i > t)
    \le
    \measure\{ Y|_{[0,1]} > t\} .
\end{equation}
This has an elementary proof --- 
see for example \cite[Proposition 2.1]{hitczenko-montgomery-smith} or
\cite{gine-zinn}.  
Thus
\begin{equation}
\label{e L_1 l_infty}
   \snormo{\snormo{(X_i)}_\infty}_1 \approx \int_0^1 Y(t) \, dt .
\end{equation}
\end{proof}

For each integer $1\le m\le n$, let $k_m$ denote the sequence space
$\snormo{(x_i)}_{k_m} = \sum_{i=1}^m x^*_i$.

\begin{lemma}
\label{l L_1 k_m}
For each positive integer $m$, 
equation~(\ref{e main}) is true if $M = L_1$ and $N = k_m$, with
constants of approximation independent of $m$.
\end{lemma}

\begin{proof}
Let $I_1$, $I_2,\dots,$\ $I_n$ be $\{0,1\}$-valued
independent random variables that
are also independent of $(X_i)$, where $\Pr(I_i = 1) = 1/m$.  
Applying equation~(\ref{e L_1 l_infty})
to the sequence $(I_i X_i)$ we obtain
\begin{equation}
\label{e Ii}
   \snormo{\snormo{(I_i X_i)}_\infty}_1 
   \approx \int_0^1 Y(mt) \, dt
   \approx \frac1m \left( \int_0^1 Y(t) \, dt + 
   \snormo{(Y(i)}_{k_m} \right).
\end{equation}
Next, let $\mathcal M$ denote the $\sigma$-field generated by $(X_i)$.
Then applying equation~(\ref{e Ii}), we see that
\[
   \E(\snormo{(I_i X_i)}_\infty | \mathcal M)
   \approx \frac1m \snormo{(X_i)}_{k_m} . 
\]
Thus we also obtain that
\[
   \snormo{\snormo{(I_i X_i)}_\infty}_1 
   = \snormo{ \E(\snormo{(I_i X_i)}_\infty | \mathcal M) }_1
   \approx
   \frac1m \snormo{\snormo{(X_i)}_{k_m}}_1 .
\]
The result follows.
\end{proof}

Let $P$ denote the space of functions $f$ on $[0,n]$ for which its quasi-norm
\[
   \snormo f_P = \snormo{f^\#|_{[0,1]}}_M + \snormo{(f^\#(i))}_N
\]
is finite.  
In fact this quasi-norm is equivalent to a norm, \emph{viz},
$\snormo f_{P'} = \snormo{f^\#|_{[0,1]}}_M 
    + \normo{\left(\int_{i-1}^i f^\#(t) \, dt\right)}_N$ (see for example
\cite[Section~7]{montgomery-smith-semenov}).  
However we will content ourselves
with proving the following statement.

\begin{lemma}
\label{l dilate P} For any function $f$ on $[0,n]$ we have
$\snormo{f(\cdot/100)}_P \le 200 \snormo f_P$.
\end{lemma}

\begin{proof}
First, since $M$ satisfies the triangle inequality, it follows
that 
$\snormo{f(\cdot/100)^\#|_{[0,1]}}_M \le 100\snormo{f^\#|_{[0,1/100]}}_M$.  
Next, since $f^\#(i/100) \le f^\#([i/100])$, where $[t]$ denotes the integer
part of $t$, we see that
\begin{eqnarray*}
   \snormo{(f^\#(i/100))}_N 
   &\le& 100 \snormo{(f^\#(i))}_N + \sum_{i=1}^{99} f^\#(i/100) \\
   &\le& 100 \snormo{(f^\#(i))}_N + 100 \int_0^1 f^\#(t) \, dt \\
   &\le& 100 \snormo f_P .
\end{eqnarray*}
\end{proof}

Finally we need to cite a couple of results.
For the case we need, $p=1$, the first result is essentially an 
immediate corollary of the 
Hoffmann-J{\o}rgensen inequality \cite{hoffmann-jorgensen}, at least
in the form found in \cite[Proposition~1.3.2]{kwapien-woyczynski}, 
and inequality~(\ref{e max in Pr}).
However we find an explicit reference to what we need in
\cite[Theorem 6.1]{hitczenko-montgomery-smith}.
The second result is
\cite[Theorem 7.1]{hitczenko-montgomery-smith}.
These concern maximal sums of vector valued random variables
$U = \max_k \normo{\sum_{i=1}^k Z_i}$, where $Z_1$, $Z_2,\dots,$\
$Z_n$ are Banach-valued independent random variables.
Let $V\colon[0,1]\to[0,\infty]$ be defined so that
\[
   \measure\{V>t\} = \min\left\{ 1 ,
   \sum_{i=1}^n \Pr(\snormo{Z_i} > t) \right\} .
\]

\begin{thm}
\label{t lp}
If
$p \ge 1$, 
then
$ \snormo U_p \approx U^\#(e^{-p}/4) + \snormo V_p $.
\end{thm}

\begin{thm}
\label{t r.i.}
Suppose that $L_q$ embeds continuously into $M$ via the natural embedding,
where $1 \le q < \infty$.  
Then
$ \snormo U_M  \approx \normo U_1 + \normo V_M $,
where the constant of approximation depends only upon $q$ and the embedding
constant.
\end{thm}

\begin{proof}[Proof of Theorem~\ref{t main}]
Let us first show the lower bound.  
Here the proof is very similar
to the proof of \cite[Theorem 27]{montgomery-smith-semenov}.
We know that
\begin{eqnarray*}
   \snormo{(x_i)}_N
   &=&
   \sup_{\snormo y_{N^*} \le 1}
   \sum_{i=1}^n x_i^* y_i^*  \\
   &=&
   \sup_{\normo y_{N^*} \le 1}
   \sum_{m=1}^n (y_m^*-y_{m+1}^*) \snormo{(x_i)}_{k_m} ,
\end{eqnarray*}
where by convention $y^*_{n+1} = 0$, and $N^*$ denotes the dual space
to $N$.  
From this, we immediately see
that
\begin{eqnarray*}
   \E \snormo{(X_i)}_N
   &\ge& 
   \sup_{\normo y_{N^*} \le 1}
   \sum_{m=1}^n (y_m^*-y_{m+1}^*)
   \E\snormo{(X_i)}_{k_m} \\
   &\approx&
   \sup_{\normo y_{N^*} \le 1}
   \sum_{m=1}^n (y_m^*-y_{m+1}^*)
   \left(
   \int_0^1 Y(t) \, dt
   +
   \snormo{(Y(i))}_{k_m}
   \right)  \\
   &\approx&
   \int_0^1 Y(t) \, dt
   +
   \snormo{(Y(i))}_N 
\end{eqnarray*}
since $y_1^* \le 1$ whenever $\normo y_{N^*} \le 1$.
To finish the lower bound, we see that
\[ 2\snormo{ \snormo{(X_i)}_N }_M
   \ge \snormo{\snormo{(X_i)}_\infty}_M + \E \snormo{(X_i)}_N , 
\]
and the result follows by equation~(\ref{e max in Pr}).

Now let us focus on the upper bound.  
Really the first part of this proof
follows by an inequality obtained independently by van Zuijlen 
\cite{van zuijlen 1}, \cite{van zuijlen 2}, \cite{van zuijlen 3},
and Marcus and Pisier \cite{marcus-pisier}.  
But we
shall provide a self contained proof that is essentially a copy of 
the proof of this same result that may be found in
\cite[Theorem 5.1]{marcus-zinn}.  
From Lemma~\ref{l dilate P},
it follows that $ \normo{Y(\cdot/100)}_P \le 200 \normo Y_P$.
We have that
\begin{eqnarray*}
   \Pr( \snormo{(X_i)}_N > 200 \normo Y_P )
   &\le&
   \Pr( \snormo{(X_i)}_N > \normo{Y(\cdot/100)}_P ) \\
   &\le&
   \Pr( \snormo{(X_i)}_N > \normo{(Y(i/100))}_N ) \\
   &\le&
   \Pr( \exists i \colon X_i^* > Y(i/100) ) \\
   &\le&
   \sum_{i=1}^n
   \Pr( X_i^* > Y(i/100) ) \\
   &\le&
   \sum_{i=1}^n
   \sum_{j_1<j_2<\cdots<j_i}
   \prod_{k=1}^i \Pr( X_{j_k} > Y(i/100) ) \\
   &\le&
   \sum_{i=1}^n
   \frac1{i!} \left(\sum_{j=1}^n \Pr(X_j > Y(i/100))\right)^i \\
   &\le&
   \sum_{i=1}^n
   \frac{i^i}{100^i i!} \\
   &\le&
   \frac1{4e} ,
\end{eqnarray*}
that is to say, $(\normo{(X_i)}_N)^\#(1/4e) \le 200 \normo Y_P$.

Now we may apply Theorems~\ref{t lp} and~\ref{t r.i.} to 
$Z_i = X_i e_i \in N$, where $e_i$ denotes the $i$th unit vector.
In that case we see that $U = \snormo{(X_i)}_N$, and $V = Y|_{[0,1]}$,
and the result follows.
\end{proof}

\section{Application to Orlicz Spaces}

In this section we will recover some of the results of
Gordon, Litvak, Sch\"utt and Werner \cite{gordon et al}.

\begin{lemma}
\label{l equiv}
Suppose that $M$ and $N$ are Orlicz spaces constructed
from Orlicz functions $\Phi$ and $\Psi$ respectively.  
Define a function
\[
   \Theta(x) = \left\{
   \begin{array}{cl}
     \Psi(x) & \text{if }0\le x \le 1\\
     \Phi(x) & \text{if }x \ge 1.
   \end{array} \right.
\]
Then $P$ is equivalent to the Orlicz space $L_\Theta$.
\end{lemma}

\begin{proof}
Note that because of the normalization on $M$ and $N$ that 
$\Phi(1) = \Psi(1) = 1$.  
Also while
$\Theta$ need not be an Orlicz function, 
it does satisfy the property that $\Theta(x)/x$ is an increasing 
function, and hence it is easily seen to
be equivalent
to the Orlicz function:
$\tilde\Theta(x) = \int_0^x \frac{\Theta(t)} t \, dt $.

Suppose that $\snormo f_{L_\Theta} \le 1$, that is
$\int_0^n \Theta(f^\#(t)) \, dt \le 1 $.
Then in particular $f^\#(1) \le \int_0^1 \Theta(f^\#(t)) \, dt \le 1 $.
Thus
\[
   \sum_{i=1}^n \Psi(f^\#(i)) \le \Theta(f^\#(1)) + 
     \int_1^n \Theta(f^\#(t)) \, dt \le 2 ,
\]
and so $\sum_{i=1}^n \Psi(f^\#(i)/2) \le 1$, that is 
$\snormo{(f^\#(i))}_N \le 2$.
Also, if $a = \measure\{f>1\}$ (so $a\le 1$), then
\[
   \int_0^1 \Phi(f^\#(t)) \, dt \le
   \int_0^a \Theta(f^\#(t)) \, dt + (1-a) \le 2,
\]
that is, $\snormo{f^\#|_{[0,1]}}_M \le 2$.  
Therefore $\snormo f_P \le 4$.

Now suppose that $\snormo f_P \le 1$.  
Again we see that $f^\#(1) \le 1$,
and $a = \measure\{f>1\} \le 1$.  
Hence
\[
   \int_0^n \Theta(f^\#(t)) \, dt \le
   \int_0^a \Phi(f^\#(t)) \, dt + (1-a) + 
   \sum_{i=1}^n \Psi(f^\#(i)) \le 3 .
\]
Since $\Theta(x/3) \le \Theta(x)/3$, it follows that
$\snormo f_{L_\Theta} \le 3$.
\end{proof}

Now we will give a formulation of one of the results of
\cite{gordon et al}, that gives a formula in the case that $M$ and
$N$ are Orlicz spaces.  
While the formula presented here may appear
superficially different than the formula given in \cite{gordon et al},
a short argument shows that it is equivalent (at least in
the case discussed in \cite{gordon et al}, that is, when $M = L_1$).

\begin{thm}
\label{t orlicz}
Suppose that $\xi_1$, 
$\xi_2,\dots,$\ $\xi_n$ are identically distributed random variables,
and that $M$ and $N$ are Orlicz spaces where 
$M$ is constructed from an Orlicz function
$\Phi$ satisfying 
\[
   \limsup_{x\to\infty} \log \Phi(x)/\log(x) < \infty .
\]
Then there exists an Orlicz
function $\Lambda$, equivalent to the function
$x\mapsto \E(\Theta(x \xi_1))$, where $\Theta$ is the function\
constructed in Lemma~\ref{l equiv}, such that for all real numbers
$a_1$, $a_2,\dots,$\ $a_n$ we have
\[
   \snormo{\snormo{(a_i \xi_i)}_N}_M \approx \snormo{(a_i)}_{L_\Lambda} .
\]
\end{thm}

\begin{proof}
The condition on $\Phi$ tells us that there exists positive constants 
$c$ and $q$ such
that
$\Phi(x) \le c x^q$ for sufficiently large $x$, that is,
$L_q$ embeds into $M$.
Thus by Theorem~\ref{t main} and Lemma~\ref{l equiv} we see that
\[
   \snormo{\snormo{(a_i \xi_i)}_N}_M \approx
   \snormo Y_\Theta .
\]
But 
\[
   \snormo Y_\Theta 
   =
   \inf\{\lambda>0 \colon \sum_{i=1}^n \E(\Theta(a_i \xi_i/\lambda)) \le 1\}
   =
   \snormo{(a_i)}_\Lambda .
\]
\end{proof}

Now we will give another proof of the following result that appears in 
\cite{gordon et al}.  
This paper also gives many other examples like this
that are interesting.

\begin{cor}
Suppose that $\gamma_1$, 
$\gamma_2,\dots,$\ $\gamma_n$ are identically distributed 
normalized Gaussian
random variables,
and let $1 \le m \le n$ be an integer.
Let $\Lambda$ be an Orlicz function equivalent to $x e^{-1/(mx)^2}$.
Then for all real numbers
$a_1$, $a_2,\dots,$\ $a_n$ we have
\begin{equation}
\label{e Gauss}
   \snormo{\snormo{(a_i \gamma_i)}_{k_m}}_1 \approx \snormo{(a_i)}_{L_\Lambda} 
   \approx \sum_{i=1}^m a^*_i 
   + m \sup_{1 \le i \le n/m} a^*_{m i} \sqrt{1+\log i} .
\end{equation}
\end{cor}

\begin{proof}
An easy argument shows that if $M=L_1$ and $N=k_m$, then
the Orlicz function $\Theta$ constructed
in Lemma~\ref{l equiv} is equivalent to the 
function $x \mapsto (x-1/m)^+$.  
Then the function
$\Lambda(x)$ constructed in Theorem~\ref{t orlicz} is equivalent to
\[
\sqrt{\frac2\pi} \int_0^\infty (x t-1/m)^+ e^{-t^2/2} \, dt
=
x \sqrt{\frac2\pi} \int_0^\infty (t-1/mx)^+ e^{-t^2/2} \, dt ,
\]
and the rest of the left approximation of (\ref{e Gauss}) 
follows by simple calculations.

To see the second approximation of (\ref{e Gauss}), note that for 
$0 \le x \le m$ that
$\Lambda(x)$ is equivalent to $x$, and that for $x>m$ that $\Lambda(x)$
is equivalent to $m e^{1-1/(mx)^2}$.  
Then 
by an argument similar to the
proof of Lemma~\ref{l equiv}, we see that 
\[ 
   \normo{(a_i)}_{L_\Lambda} \approx
   \normo{(a^*_i)_{1 \le i \le m}}_1
   +
   m \normo{(a^*_{m i})}_{L_{\exp(1-1/x^2)}} .
\]
Finally we need to show that
\[
   \normo{(b^*_{i})}_{L_{\exp(1-1/x^2)}} 
   \approx
   \sup_{i} b^*_{i} \sqrt{1+\log i} .
\]
This is essentially a sequential version of results from
\cite{bennett-rudnick}.  
Suppose that 
$\normo{(b^*_{i})}_{L_{\exp(1-1/x^2)}} \le 1$.  
Then for any positive integer 
$i$,
\[
   i \exp(1-1/(b^*_i)^2) \le \sum_j \exp(1-1/(b^*_j)^2) \le 1,
\]
that is, $b^*_{i} \sqrt{1+\log i} \le 1$.  
Conversely, 
if $\sup_i b^*_i \sqrt{1+\log i} \le 1$,
then
\[
   \sum_i \exp(1-1/(b^*_i/2)^2)
   \le
   \sum_i \exp(-3-4\log i)
   = 
   e^{-3} \sum_i \frac1{i^4} \le 1,
\]
that is, $\normo{(b^*_{i})}_{L_{\exp(1-1/x^2)}} \le 2$.
\end{proof}

We remark that a similar argument shows that
if $\xi_1$, $\xi_2,\dots,$\ $\xi_n$ are identically distributed
random variables with $E|\xi_1| = 1$, then
\[
   \snormo{\snormo{(a_i \xi_i)}_{k_m}}_1 
   \approx \sum_{i=1}^m a^*_i  + 
    m \snormo{\snormo{(a^*_{m i} \xi_i)}_\infty}_1 .
\]
\bigskip

Finally let us finish with another remark.  
In \cite{gordon et al}, the authors
showed in the case that $M = L_1$
that their upper bound held even if the random variables were
not independent.  
This can also hold in our more general case.
In \cite{montgomery-smith-semenov} was introduced the concept of what
it means for a rearrangement invariant space to be $D^*$-convex.
This property is held, for example, by all Orlicz spaces.
Following the proof of \cite[Theorem 27]{montgomery-smith-semenov},
it can be shown that equation~(\ref{e main}) holds even if the
sequence $(X_i)$ is not necessarily independent, as long as 
$M = L_1$ and $P$ is $D^*$-convex.  
It is easy to see from the definition that the condition that
$P$ be $D^*$-convex cannot be dropped.
We leave the details to the
interested reader.

\end{document}